\documentclass[1p,times,sort]{elsarticle}
\usepackage[utf8]{inputenc}
\usepackage{amsmath,amsfonts,amssymb,amsthm}
\usepackage[all]{xy}
\usepackage{graphicx,subfig,aliascnt}
\usepackage[colorlinks=true]{hyperref}
\usepackage[ruled,linesnumbered,noend]{algorithm2e}

\journal{arXiv}

\theoremstyle{plain}
\newtheorem{theorem}{Theorem}[section]
\newtheorem{proposition}[theorem]{Proposition}
\newtheorem{corollary}[theorem]{Corollary}
\newtheorem{lemma}[theorem]{Lemma}
\theoremstyle{definition}
\newtheorem{definition}[theorem]{Definition}
\newtheorem{example}[theorem]{Example}
\newtheorem{remark}[theorem]{Remark}
\newtheorem{notation}[theorem]{Notation}

\theoremstyle{plain}
\newtheorem*{theorem*}{Theorem}
\newtheorem*{proposition*}{Proposition}
\newtheorem*{corollary*}{Corollary}
\newtheorem*{lemma*}{Lemma}
\theoremstyle{definition}
\newtheorem*{definition*}{Definition}
\newtheorem*{example*}{Example}
\newtheorem*{remark*}{Remark}
\newtheorem*{notation*}{Notation}

\begin{document}

\begin{frontmatter}

\title{Solving congruence equations using Bernstein forms}

\author[dm]{C\'esar Massri\corref{correspondencia}\fnref{financiado}}
\ead{cmassri@dm.uba.ar}
\address[dm]{Department of Mathematics, University CAECE, Buenos Aires, Argentina}
\cortext[correspondencia]{Address for correspondence: Department of Mathematics, FCEN,
University of Buenos Aires, Argentina. Postal address: 1428. Phone number: 54-11-4576-3335.}
\fntext[financiado]{The author was fully supported by CONICET, IMAS, Buenos Aires, Argentina}

\author[dc]{Manuel Dubinsky\fnref{financiado2}}
\ead{mdubinsky@undav.edu.ar}
\address[dc]{Computer Engineering Department, University of Avellaneda, Argentina}
\fntext[financiado2]{The author was fully supported by UndavCyT 2014 project.}

\begin{abstract}
We present a subdivision method to solve systems of congruence
equations. 
This method is inspired in a subdivision method,
based on Bernstein forms, 
to solve systems of polynomial inequalities in several variables and arbitrary
degrees. The proposed method is exponential in the number of variables.
\end{abstract}

\end{frontmatter}

\section*{Introduction}

\subsection*{Overview of the problem}

An important aspect within the study of algebraic varieties
over characteristic $p$ is to know the number of points that a variety has.
This led to the statement of Weil's conjectures
proved by Deligne thanks to the analysis made by Grothendieck and others.
In this article we are interested in computing, not only
the number of points in an algebraic variety modulo $p$
but also, the coordinates of these points. We propose a computational 
method to find solutions to the following
system of congruence equations

\[
\left\{\begin{array}{rcl}
h_1(x_1,\ldots,x_n)&\equiv 0&\mbox{mod }(p_1)\\
&\vdots\\
h_r(x_1,\ldots,x_n)&\equiv 0&\mbox{mod }(p_r)
\end{array}\right .
\]
where $h_1,\ldots,h_r\in\mathbb{Z}[x_1,\ldots,x_n]$ and 
$p_1,\ldots,p_r\in\mathbb{Z}_{\ge 2}$. To achieve this goal we 
first analyze a known method in nonlinear programming and then adapt 
it to our needs.

Polynomial programming is the process of solving
a system of equalities and inequalities between
polynomial functions along with an objective function
to be maximized (or minimized).
Typically, such a system can be described as:
\begin{align*}
\text{maximize }& q_1(x)\\
\text{subject to }& q_i(x)\ge 0,\quad 2\le i\le r.
\end{align*}
The polynomials $q_i$ are assumed to have real coefficients
but in this article
we study integer polynomial programming
assuming $q_1,\ldots,q_r\in\mathbb{Q}[x_1,\ldots,x_n]$ and $x\in\mathbb{Z}^n$.

Since the resolution of Hilbert's Tenth problem
by Matiyasevich, \citep{MR1244324}, several
results appeared proving that 
there cannot exist any general algorithm to solve
integer polynomial programming (it is non computable).
For example, in \citep{jeroslow}, the author
proves that the problem of minimizing a linear
form over quadratic constraints in
integer variables is not computable by a recursive function.
Hence, in order to avoid this phenomena,
it is necessary to add some constrains to the original problem.
Given that we are interested in applications to
congruence equations we bound the variables to some region
$D=[a_1,b_1]\times\ldots\times[a_n,b_n]\subseteq\mathbb{R}^n$
restricting the problem of adapting an integer polynomial
programming
\begin{align*}
\text{maximize }& q_1(x)\\
\text{subject to }& q_i(x)\ge 0,\quad 2\le i\le r\\
\text{and }&x\in D\cap\mathbb{Z}^n.
\end{align*}
to a solver of congruence equations
\[
\left\{\begin{array}{rcl}
h_1(x)&\equiv 0&\mbox{mod }(p_1)\\
&\vdots\\
h_r(x)&\equiv 0&\mbox{mod }(p_r)
\end{array}\right .
\]
where $x\in D\cap \mathbb{Z}^n$, $D=[0,k]^n$ and
$k$ is the least common multiple of $p_1,\ldots,p_r$,
$k=\mbox{lcm}(p_1,\ldots,p_r)$.

\subsection*{Existing work}

There are many local and global methods to solve a Diophantine equation.
The first thing to do is to check whether the problem has a solution locally.
If it has none, the problem is solved since our equation has
no solution, but if it has, the local information that we obtain
may already completely solve the problem, or may not
so it is necessary to study the equation
globally. The global methods include
sophisticated tools from number theory, for example
factorization over a number field, computations of class groups and unit groups,
Diophantine approximation techniques, 
linear forms in logarithms, elliptic curves,
modular forms and Galois representations and more.

Regarding congruence equations we can say that
there are general methods for solving both linear and quadratic congruence equations
however solutions to the general polynomial congruence is intractable.
Standard algorithms can be found at \citep{MR1689189}
and also at \citep{MR2312337,MR2312338}. Several software packages
implement these techniques, see for example MagmaSoft, SageMath, Kant/Kash.

\

Integer polynomial
programming is a subclass of
the more general class \emph{Mixed integer
nonlinear programming} (MINLP).
Approaches to solve these problems may be classified
as either stochastic or deterministic.
Stochastic techniques employ some element of
randomness in their search for the global optimum
and consequently rely on a statistical argument to
prove convergence to the global optimum, see \cite{MR1263591}.
Deterministic methods have the advantage that they
provide a rigorous guarantee of global optimality of
any solution produced within a specified tolerance.
Examples of deterministic techniques are cutting plane
methods, interval methods, primal-dual methods,
outer approximation, generating functions, etc.
For instance using generating functions
of convex sets, it is proven in \citep{MR1992831}
that there exists a polynomial-time algorithm
to find the maximum value of the objective
function without knowing where this maximum is achieved.
Also, a very successful methodology
is the spatial branch and bound algorithm \citep{MR1668773,MR2640549}.
The aim of branch and bound techniques,
is to divide the problem into subclasses to be solved with
convex or linear approximations, which are
powerful and flexible for
modeling and solving decision problems.
Despite their widespread use, few available
software packages provide any guarantee of correct answers or
certification of results.
Possible inaccuracy is caused by the use of floating-point numbers.
Floating-point calculations require the use of built-in tolerances
for testing feasibility and optimality,
and can lead to calculation errors in the solution of linear programming
relaxations and in the methods used for creating cutting planes
to improve these relaxations.
See in \citep{MR3102373} where the author gave a method
to guarantee the result avoiding numerical errors.

There is a large variety of commercial
and academic noncommercial software, but
only a few of those software packages solve general
non-convex problems to global optimality.
For a review see \citep{bussieck2010minlp},
\citep{MR3473492}, \citep[\S 1]{book:747990} and/or \citep{MR2869850}.
It is worth mention that there is a large literature reviewing
different methods (and software) developed to solve
convex and non-convex MINLP problems, furthermore,
every two years, Robert Fourer publishes a list of currently
available codes in the field of linear
and integer programming, the latest edition can
be found at the web page of OR/MS-Today, Software Survey
section, \citep{fourer}.

In this article we use Bernstein forms. The advantage of using Bernstein 
forms is that they produce robust and reliable algorithms. 
They are used to solve systems of polynomial equations and also integer 
polynomial programming,
see for example \citep{garloff1993bernstein,nataraj2007new,patil2014improved}.
We review an algorithm using Bernstein forms to solve 
integer polynomial programming and then we adapt it to solve congruence 
problems.

\

Before recalling the next theorem, let 
us introduce some notations.
Let $(d_1,\ldots,d_n)$ be a multi-degree, $s=(d_1+1)\ldots(d_n+1)$ and
let $D=[a_1,a_1+2^{k_1}]\times\ldots\times[a_n,a_n+2^{k_n}]$,
$(a_1,\ldots,a_n)\in\mathbb{Z}^n$,
$(k_1,\ldots,k_n)\in\mathbb{N}_0^n$.
Choose some index $j$, $1\le j\le n$, and then divide
into two halves the $j$-side of $D$,
\begin{align*}
D_j^L=&\prod_{i=1}^{j-1}[a_i,a_i+2^{k_i}]\times[a_j,a_j+2^{k_j-1}]\times\prod_{i=j+1}^{n}[a_i,a_i+2^{k_i}],\\
D_j^R=&\prod_{i=1}^{j-1}[a_i,a_i+2^{k_i}]\times[a_j+2^{k_j-1},a_j+2^{k_j}]\times\prod_{i=j+1}^{n}[a_i,a_i+2^{k_i}].
\end{align*}
\begin{theorem*}
For every $q_1,\ldots,q_r\in
\mathbb{Z}[x_1,\ldots,x_n]_{\le(d_1,\ldots,d_n)}$
there exists a rectangular matrix $v\in\mathbb{Q}^{s\times r}$
such that
the values of $q_i$ over $D_j^L$ (resp. $D_j^R$) are bounded by
the minimum and the maximum of
$\{w^L_{ki}\}_{k=1}^s$ (resp. $\{w^R_{ki}\}_{k=1}^s$),
where $w^L=M_j^Lv$, $w^R=M_j^Rv$, $1\le i\le r$. The matrices
$\{M_j^L,M_j^R\}_{j=1}^n$ depend on $(d_1,\ldots,d_n)$.

Also, the value of $q_i$ at $(a_1,\ldots,a_n)$ is equal to $w^L_{1i}$
and the value of $q_i$ at $(a_1,\ldots,a_{j-1},a_j+2^{k_j-1},a_{j+1},\ldots,a_n)$
is equal to $w^R_{1i}$, $1\le i\le r$.
\end{theorem*}

This theorem is related to the resolution of integer polynomial 
programming as follows.
The subdivision process tests if a given
region may have solutions or definitely not.
If it may, then the region is divided and the process starts again, 
otherwise it is rejected. At the end, the process
produces candidates of the form $\prod_{i=1}^n[z_i,z_i+1]$
and test if $(z_1,\ldots,z_n)$ is a solution
of the system. Both stages of this algorithm depend
on the previous theorem.

\subsection*{Main result}

In this article we present a subdivision
method based on Bernstein forms
to compute the solutions of a system of congruence equations.
The idea of using Bernstein forms
to make a ``solver'' is not new, see for example
\citep{garloff1993bernstein,nataraj2007new,patil2014improved} and
\citep{MR3079719,MR2499511,Elber:2001:GCS:376957.376958}.
One of our contribution is
the application of Bernstein forms
to congruence equations. The key ingredient is that we managed to make
the subdivision process and the process of rewriting
the system in Bernstein form as a matrix multiplication.

An advantage of Bernstein forms is that the algorithm can test
the existence of a solution in a large region by
testing bounds on the coordinates of a matrix.
An important feature of the algorithm
is that it does not rely on numerical computations
and as a result, it gives a certified answer,
see Theorem \ref{answer2}.
The algorithm works by
performing a matrix multiplication at each iteration
with triangular rational matrices $\{M_i^L,M_i^R\}_{i=1}^n$.
These matrices
depend on the multi-degree of the system $(d_1,\ldots,d_n)$
and we prove that it is convenient to have these $2n$
matrices previously computed
for a large multi-degree $(d_1,\ldots,d_n)$,
see Theorem \ref{comple2}. It is worth mention
that the algorithm can work with any system of congruence equations, 
independently if the numbers $p_1,\ldots,p_r$ are primes or not.

\

The first contribution that we present is
Theorem \ref{comple} 
where we give the expected number
of matrix multiplications that the algorithm 
for integer polynomial programming requires.
It is easily seen that this algorithm has an
exponential complexity, but in practice,
we prove that the complexity depends on
the ``size'' of the real solution of
the system, see Proposition \ref{elev2}.
This ``size'' can be estimated
using tools from algebraic geometry (dimension, degree)
and in particular, Theorem \ref{comple} implies
that the algorithm is faster
if we add more equations, that is, if there are less solutions
to the system.
We used concepts and
results from the theory of Branching processes
to prove Theorem \ref{comple},
\begin{theorem*}
Let $\lambda$ be the complexity number (Definition \ref{complex-number}) of
$\{q_1,\ldots,q_r\}$ over $D=\prod_{i=1}^n[a_i,a_i+2^{k_i}]$
and let $K=k_1+\ldots+k_n$.
Then, the expected number of matrix multiplications produced by the algorithm is
\[
\sum_{i=0}^K
2^{i}\prod_{j=0}^{i-1}\left(1-(1-\lambda)^{2^{K-j}}\right)-1.
\]
\end{theorem*}

Regarding systems of congruence equations we first give
Lemma \ref{lemma-hasSol} that characterizes existence of solutions in a real interval.
Then, we give the pseudocode of {\tt SolveCE} and prove
Theorem \ref{answer2}

\begin{theorem*}
The output of {\tt SolveCE}($(p_1,\ldots,p_r)$,$v$, $(k,\ldots,k)$,$(0,\ldots,0)$)
is the set
\[
\left\{x\in D\cap\mathbb{Z}^{n}\,\colon\,
p_1|h_1(x),\ldots,
p_r|h_r(x)\right\},
\]
where
$D=[0,2^{k}]\times\ldots\times[0,2^{k}]$,
$k=\lceil\log_2(\mbox{lcm}(p_1,\ldots,p_r))\rceil$ and
$v$ is the matrix associated to $h_1,\ldots,h_r$.
\end{theorem*}

Finally, we compute the expected number of matrix multiplications
produced by {\tt SolveCE}. This number depends on a parameter 
$\lambda$ (similar to the complexity number),

\begin{theorem*}
Let $h_1,\ldots,h_r\in\mathbb{Z}[x_1,\ldots,x_n]$,
$p_1,\ldots,p_r\in\mathbb{Z}_{\ge 2}$, 
$k=\lceil\log_2(\mbox{lcm}(p_1,\ldots,p_r))\rceil$.
Then, the expected number of matrix multiplications produced by {\tt SolveCE}
in $D=\prod_{i=1}^n[0,2^{k}]$ is
\[
\sum_{i=0}^{nk}
2^{i}\prod_{j=0}^{i-1}\left(1-(1-\lambda)^{2^{nk-j}}\right)-1.
\]
This number is bounded between
a polynomial and an exponential expression in $\lambda$.
\end{theorem*}

\subsection*{Summary}

In \textbf{Section 1} we give some preliminary
results; we recall the concept of Bernstein Bases
and degree elevation, fixing a multi-degree $(d_1,\ldots,d_n)$
we associate to every region $D$
a square matrix and to every system $\{q_1,\ldots,q_r\}$
a rectangular matrix. Finally, with these results, we prove
in Theorem \ref{main} our main tool.
In \textbf{Section 2} we present the algorithm
to solve integer polynomial programming {\tt SolveIPP}
and in \textbf{Section 3} we study its expected complexity
and give several examples computed 
using our implementation in SageMath \citep{sage}.
Finally in \textbf{Section 4} we present the
pseudocode {\tt SolveCE} to solve congruence equations
and study its expected complexity. We compare
in a table {\tt SolveCE} with the brute-force algorithm (Definition \ref{bf}).

\section{Preliminaries}

In this and the next section we present a known algorithm to solve
integer polynomial programming. Later we adapt it
to solve congruence equations. 

We are interested in the integral points of
semi-algebraic sets defined by
a system of multivariate rational polynomial equations
in $D=[a_1,b_1]\times\ldots\times[a_n,b_n]\subseteq\mathbb{R}^n$.
Specifically,
let $q_1\in\mathbb{Q}[x_1,\ldots,x_n]$ be some multivariate
polynomial and let $\mathcal{S'}$ be the intersection
between $D\cap\mathbb{Z}^n$ and some semi-algebraic set defined by
some polynomials $q_2,\ldots,q_r\in \mathbb{Q}[x_1,\ldots,x_n]$. We
want to find the maximum $\theta$ of $q_1$ over $\mathcal{S'}$
and more general, to describe the set $\mathcal{S}$,
\[
\mathcal{S}=\{\mathbf{z}\in\mathcal{S'}\,\colon\,q_1(\mathbf{z})=\theta\}.
\]
\begin{notation}
First of all, we want to give a unifying description of the set
$\mathcal{S'}$.
It is defined
as the intersection of equalities and inequalities.
The first simplification to make is to assume that the polynomials
have integer coefficients and then,
since we are working over the integers, we can
always assume that $\mathcal{S'}$
is defined only with inequalities of the form $q\ge 0$.
Concretely, we can make the following changes
for $q,q'\in\mathbb{Z}[x_1,\ldots,x_n]$:
\begin{itemize}
\item Replace the inequality $q\le q'$ with the inequality $q'-q\ge 0$.
\item Replace the equality $q=q'$ with the inequalities $q-q'\ge 0$ and $q'-q\ge0$.
\item Replace the inequality $q>q'$ with the inequality $q-q'-1\ge 0$.
\item Replace the inequality $q<q'$ with the inequality $q'-q-1\ge 0$.
\end{itemize}
Then, without loss of generality,
we can assume that $\mathcal{S'}$ is given as,
\[
\mathcal{S'}=\left\{\mathbf{z}\in D\cap\mathbb{Z}^n\,\colon\,
q_2(\mathbf{z})\ge 0,\,\ldots\,,q_r(\mathbf{z})\ge 0\right\},
\]
where $q_2,\ldots,q_r\in\mathbb{Z}[x_1,\ldots,x_n]$ are
some multivariate polynomials.
Also, we can always assume that $q_1$ has integer coefficients
and that we are searching for the maximum of $q_1$ over
$\mathcal{S'}$ making the following changes,
\begin{itemize}
\item Replace $\min(q_1)$ with $\max(-q_1)$.
\item If there is no $q_1$ to maximize, then maximize $q_1\equiv 0$.
\end{itemize}
Then, the solutions of the system is the set $\mathcal{S}$,
\[
\mathcal{S}=\left\{\mathbf{z}\in \mathcal{S'}\,\colon
q_1(\mathbf{z})=\max_{\mathcal{S'}}(q_1)\right\}.
\]
\qed
\end{notation}

Now that we have a unifying way to present the system,
let us define (or recall) Bernstein bases.
To simplify the notation, we use standard multi-index notation.
Letters in boldface represent multi-indexes. For example
$\mathbf{d}=(d_1,\ldots,d_n)$.
\begin{definition}[Bernstein Basis]
According to \cite[Proposition 2]{MR3079719},
any polynomial $p\in\mathbb{R}[x_1,\ldots,x_n]$
of multi-degree $(d_1,\ldots,d_n)$ can be written
in Bernstein form,
\[
p(\mathbf{x})=
p(x_1,\ldots,x_n)=\sum_{k_1=0}^{d_1}\ldots\sum_{k_n=0}^{d_n}
\beta_{k_1\ldots k_n} B_{d_1\ldots d_n,k_1\ldots k_n}(x_1,\ldots,x_n)=
\sum_{\mathbf{k}=\mathbf{0}}^\mathbf{d}
\beta_\mathbf{k}B_{\mathbf{d},\mathbf{k}}(\mathbf{x}),
\]
where $\beta_\mathbf{k}=\beta_{k_1\ldots k_n}\in\mathbb{R}$ and
\[
B_{\mathbf{d},\mathbf{k}}(\mathbf{x})=
\binom{d_1}{k_1} x_1^{k_1}(1-x_1)^{d_1-k_1}\ldots
\binom{d_n}{k_n} x_n^{k_n}(1-x_n)^{d_n-k_n}=
\binom{\mathbf{d}}{\mathbf{k}}
\mathbf{x}^\mathbf{k}(\mathbf{1}-\mathbf{x})^{\mathbf{d}-\mathbf{k}}.
\]
The polynomials $B_{\mathbf{d},\mathbf{k}}$
have rational coefficients and multi-degree
$\mathbf{d}$ for every $\mathbf{k}$.
The ordered set
\[
\{B_{\mathbf{d},\mathbf{k}}\,\colon\, \mathbf{0}\le \mathbf{k}\le\mathbf{d}\}
\]
is called
the \emph{Bernstein basis} of $\mathbb{R}[\mathbf{x}]_{\le\mathbf{d}}$
or the \emph{$\mathbf{d}$-Bernstein basis} to emphasize the multi-degree.
The order is the lexicographical order.
The numbers $\beta_\mathbf{k}$ are called \emph{Bernstein coefficients}.
\qed
\end{definition}

\begin{remark}
If a polynomial is given in monomial form
$p(\mathbf{x})=\sum_{\mathbf{i}=\mathbf{0}}^\mathbf{d} c_\mathbf{i} \mathbf{x}^\mathbf{i}$, it is easy to make a change
of basis from monomial basis to Bernstein basis.
Specifically, the Bernstein coefficients are computed as follows:
\[
\beta_\mathbf{k}=\sum_{\mathbf{i}=\mathbf{0}}^\mathbf{k} c_\mathbf{i}\frac{\binom{\mathbf{k}}{\mathbf{i}}}{\binom{\mathbf{d}}{\mathbf{i}}}
\]
Note that $\{\beta_\mathbf{k}\}\subseteq\mathbb{Q}$
if $\{c_\mathbf{i}\}\subseteq\mathbb{Q}$.
Recall from \cite[Corollary 1]{MR3079719} that
if $m_1$ (resp. $m_2$) is the minimum (resp. maximum) of the
Bernstein coefficients $\beta_\mathbf{k}$, then
we have the \emph{fundamental inequality}
\[
m_1\le p(\mathbf{x})\le m_2,
\quad\forall\mathbf{x}\in [0,1]^n.
\]
\qed
\end{remark}

\begin{lemma}[Degree elevation]\label{elev}
Let $p\in\mathbb{Z}[\mathbf{x}]_{\le\mathbf{d}}$
be a polynomial defined over $[0,1]^n$ and let $m_1,m_2\in\mathbb{R}$ be
such that $m_1<p(\mathbf{x})<m_2$
for all $\mathbf{x}\in [0,1]^n$.
Then, there exists $\mathbf{d'}$ such that for every
$\mathbf{d''}\ge\mathbf{d'}$
the Bernstein coefficients $\beta_\mathbf{k}$ of $p$
in $\mathbf{d''}$-Bernstein basis satisfy,
\[
m_1<\beta_\mathbf{k}<m_2,
\quad\forall
\mathbf{k}\le \mathbf{d''}.
\]
\end{lemma}
\begin{proof}
Given that $m_1<p(\mathbf{x})<m_2$, there exists
$\delta>0$ such that
$m_1+\delta<p(\mathbf{x})<m_2-\delta$
for all $\mathbf{x}\in [0,1]^n$.
Also, from \citep[Eq. 3.1]{MR1266028} there exists
$m\in\mathbb{N}$ such that
$|\beta_\mathbf{k}-p(\mathbf{k}/m)|\le\delta$
for all $\mathbf{k}\le \mathbf{d'}:=(m,\ldots,m)$.
Then,
\[
\beta_\mathbf{k}=\left(\beta_\mathbf{k}-p(\mathbf{k}/m)\right)+p(\mathbf{k}/m)
<\delta+(m_2-\delta)=m_2
\]
and
\[
\beta_\mathbf{k}=\left(\beta_\mathbf{k}-p(\mathbf{k}/m)\right)+p(\mathbf{k}/m)
>-\delta+(m_1+\delta)=m_1.
\]
\end{proof}

\begin{example}
Consider the following polynomial over $[0,1]$,
\[
p(x)=-\left(x-\frac{1}{2}\right)^2-\frac{1}{10}.
\]
Then $\max_{[0,1]}(p)=-1/10$.
The Bernstein coefficients in degree $2$ of $p$ are
$(-7/20,3/20,-7/20)$. Taking the
Bernstein coefficients in degree $d'=3$ we get
$(-7/20,-1/60,-1/60,-7/20)$.
In particular, for every $d''\ge 3$,
\[
-\frac{1}{10}\le \max_{k\le d''}(\beta_k)\le -\frac{1}{60}.
\]
\qed
\end{example}

Let us start analyzing the regions $D$. Our goal is to
define subdivision matrices for $[0,1]^n$ and then,
to treat any $D$ in a uniform way.

\begin{definition}[Matrix associated to a box and to a multi-degree]\label{box}
We say that $D\subseteq\mathbb{R}^n$ is a \emph{box} if there exist
$\mathbf{a},\mathbf{b}\in\mathbb{R}^n$ such that
\[
D=[a_1,b_1]\times\ldots\times[a_n,b_n].
\]
Let $\varphi:\mathbb{R}^n\rightarrow\mathbb{R}^n$
be $\varphi(\mathbf{x})=(\mathbf{b}-\mathbf{a})\mathbf{x}+\mathbf{a}$.
Clearly $\varphi$ maps bijectively $[0,1]^n$ to $D$
and defines, via
pull-back, a linear map $\mathbb{R}[\mathbf{x}]\rightarrow
\mathbb{R}[\mathbf{x}]$ given by
$q\mapsto \varphi^{*}(q):=q((\mathbf{b}-\mathbf{a})\mathbf{x}+\mathbf{a})$.

Fix a multi-degree $\mathbf{d}$.
We say that the matrix $M$ of the map $\varphi^{*}$ in
the $\mathbf{d}$-Bernstein basis
is the \emph{matrix associated to the box $D$
and to the multi-degree $\mathbf{d}$}.
The matrix $M$ is a $s\times s$-square matrix,
where $s=(d_1+1)\ldots(d_n+1)$ and
also, if $\{a_i,b_i\}_{i=1}^n\subseteq\mathbb{Q}$,
then $M\in\mathbb{Q}^{s\times s}$.
\qed
\end{definition}

\begin{notation}\label{def-m}
Fix a multi-degree $\mathbf{d}$.
Let us denote $M_{i}^L$ or $M_{i}^L(\mathbf{d})$ (resp. $M_{i}^R$ or $M_{i}^R(\mathbf{d})$) to the matrix
associated to the box
$[0,1]^{i-1}\times [0,1/2]\times[0,1]^{n-i}$
(resp. $[0,1]^{i-1}\times [1/2,1] \times[0,1]^{n-i}$)
for
$1\le i\le n$ and to $\mathbf{d}$
(see Definition \ref{box}).
Recall that if
$s=(d_1+1)\ldots(d_n+1)$, then
$M_i^L,M_i^R\in \mathbb{Q}^{s\times s}$.
\qed
\end{notation}

\begin{proposition}
The matrix $M_i^L(\mathbf{d})$ is lower triangular
and $M_i^R(\mathbf{d})$ is upper triangular,
$1\le i\le n$.
If $i\neq j$, $M_i^L$ commutes with $M_j^L$ and
with $M_j^R$, but $M_i^LM_i^R\neq M_i^RM_i^L$.
\end{proposition}
\begin{proof}
The shape of the matrices $M_i^L$ and $M_i^R$ follows from
\cite[Proposition 5]{MR3079719}.
Let us prove the commutativity.

Let $\varphi_i^L$ (resp. $\varphi_i^R$) be the unique
affine isomorphism sending
$[0,1]^n$ to $[0,1]^{i-1}\times [0,1/2]\times[0,1]^{n-i}$
(resp. to $[0,1]^{i-1}\times [1/2,1]\times[0,1]^{n-i}$).
The matrix representation of $\varphi_{i}^L$ in $\mathbf{d}$-Bernstein
basis is $M_{i}^L$ (same for $\varphi_{i}^R$ and $M_{i}^R$).
Then, we need to prove that
if $i\neq j$, then $\varphi_i^L$ commutes with $\varphi_j^L$ and
with $\varphi_j^R$, but $\varphi_i^L\varphi_i^R\neq \varphi_i^R\varphi_i^L$. Clearly $\varphi_i^L\varphi_j^L$ and
$\varphi_j^L\varphi_i^L$ (resp. $\varphi_i^L\varphi_j^R$
and $\varphi_j^R\varphi_i^L$)
are affine isomorphisms
sending $[0,1]^n$ to the same box. By uniqueness
$\varphi_i^L\varphi_j^L=\varphi_j^L\varphi_i^L$
(resp. $\varphi_i^L\varphi_j^R=\varphi_j^R\varphi_i^L$).

Given that the box associated to
$\varphi_i^L\varphi_i^R$ is
different from the box associated to
$\varphi_i^R\varphi_i^L$, we obtain
$\varphi_i^L\varphi_i^R\neq \varphi_i^R\varphi_i^L$,
\[
[0,1]^{i-1}\times [1/2,3/4]\times[0,1]^{n-i}\quad\neq\quad
[0,1]^{i-1}\times [1/4,1/2]\times[0,1]^{n-i}.
\]
\end{proof}

\begin{proposition}\label{subdiv}
Fix a multi-degree $\mathbf{d}$.
Let $k_i,l_i\in\mathbb{N}_0$ be
such that $l_i<2^{k_i}$, $1\le i\le n$.
Let $D$ be the box defined as
\[
D=\left[\frac{l_1}{2^{k_1}},\frac{l_1+1}{2^{k_1}}\right]\times\ldots\times\left[\frac{l_n}{2^{k_n}},\frac{l_n+1}{2^{k_n}}\right]
\]
and let $M$ be the matrix associated to $D$ and $\mathbf{d}$.
Then, there exists a factorization of $M$ into a product
of matrices $M_{i}^L(\mathbf{d})$ and $M_{i}^R(\mathbf{d})$.
The sum of the multiplicities of the
matrices $M_{i}^L(\mathbf{d})$ and $M_{i}^R(\mathbf{d})$
in this product is $k_i$, $1\le i\le n$.
\end{proposition}
\begin{proof}
Without loss of generality, we may assume that $n=1$ and
$D=[l/2^k,(l+1)/2^k]$. Let us call $M_1^L=L$ and $M_1^R=R$.

Induction in $k$. If $k=1$, then $l=0$ or $l=1$. Hence, $M=L$ or $M=R$.
Assume now the result for $k-1$ and let us prove it for $k$.
We have two possibilities, $l$ is even or $l+1$ is even.
If $l=2l'$ is even, then $D$ is included in $D'=[l/2^k,(l+2)/2^k]=[l'/2^{k-1},(l'+1)/2^{k-1}]$. By the inductive hypothesis, $D'$
has associated a matrix $M'$ that can be factorized
as $k-1$ products of the matrices $L$ and $R$,
$M'=L^{i_1}R^{i_2}\ldots $. Given that $D$ is equal to
the first half of $D'$, $M=LM'$. The case in which $l+1$ is even is similar.
\end{proof}

\begin{definition}
Let $F(k)$ be the set
of functions from $\{1,\ldots,k\}$ to $\{L,R\}$.
The cardinal of $F(k)$ is $2^{k}$.

Fix a multi-degree $\mathbf{d}$ and $\mathbf{k}\in\mathbb{N}_0^n$.
The \emph{set of $(\mathbf{d},\mathbf{k})$-subdivision matrices}
is defined as
\[
\left\{ M_1^{f_1(1)}\ldots M_1^{f_1(k_1)}
M_2^{f_2(1)}\ldots M_2^{f_2(k_2)}\ldots
M_n^{f_n(1)}\ldots M_n^{f_n(k_n)}
\,\colon\, (f_1,\ldots,f_n)\in F(k_1)\times\ldots\times F(k_n)
\right\}.
\]
If $k=0$, we use the convention $F(k)=\{0\}$ and $M_i^0$ is the identity matrix.
According to Proposition \ref{subdiv}, the cardinal
of this set is equal to the number of boxes obtained
by subdividing $k_i$ times in direction $x_i$ the box $[0,1]^n$.
Hence, the cardinal is equal to $2^{k_1+\ldots+k_n}$.
\qed
\end{definition}

Now that we know how to treat different boxes $D$
as matrices, let us give some results on how
to treat polynomials as vectors.

\begin{proposition}\label{mat-vec}
Let $\mathbf{d}\in\mathbb{N}_0^n$,
$s=(d_1+1)\ldots(d_n+1)$,
$D=[a_1,b_1]\times\ldots\times[a_n,b_n]$
and let $q\in\mathbb{R}[\mathbf{x}]_{\le\mathbf{d}}$
be a multivariate polynomial of multi-degree $\le\mathbf{d}$.
Then, there exist a vector $v=v(q)\in\mathbb{R}^s$ such
that,
\begin{itemize}
\item The first coefficient $v_1$ is equal to $q(a_1,\ldots,a_n)$,
\[
v_1=q(a_1,\ldots,a_n).
\]
\item The values of $q$ over $D$ are bounded by $\min(v)$ and
$\max(v)$,
\[
\min(\{v_i\}_{i=1}^s)\le q(\mathbf{x})\le \max(\{v_i\}_{i=1}^s),\quad\forall\mathbf{x}\in D.
\]
\end{itemize}
\end{proposition}
\begin{proof}
Let $\varphi:\mathbb{R}^n\rightarrow\mathbb{R}^n$
be the unique affine isomorphism
sending $[0,1]^n$ to $D$,
$\varphi(\mathbf{x})=(\mathbf{b}-\mathbf{a})\mathbf{x}+\mathbf{a}$.
Then, the pull-back of $q$ over $D$, defines a polynomial
$\varphi^*(q)=q((\mathbf{b}-\mathbf{a})\mathbf{x}+\mathbf{a})$
over $[0,1]^n$.
Let $v$ be the vector representation of
$\varphi^{*}(q)$ in $\mathbf{d}$-Bernstein basis. Then,
\[
\min(v)\le q(\mathbf{x})\le \max(v),\quad\forall\mathbf{x}\in D
\]
and from \cite[Proposition 4]{MR3079719},
\[
v_1=b_{\mathbf{0}}=\varphi^*(q)(\mathbf{0})=
q(\varphi(\mathbf{0}))=q(a_1,\ldots,a_n).
\]
\end{proof}

\begin{corollary}
Let $\mathbf{d}\in\mathbb{N}_0^n$,
$s=(d_1+1)\ldots(d_n+1)$,
$D=[a_1,b_1]\times\ldots\times[a_n,b_n]$
and let $q_1,\ldots,q_r\in\mathbb{R}[\mathbf{x}]_{\le\mathbf{d}}$
be polynomials of multi-degree $\le\mathbf{d}$.
Then, there exist a
rectangular matrix $v\in\mathbb{R}^{s\times r}$ with columns
$v_i=v_i(q_i)$ such that,
\begin{itemize}
\item The coefficient $v_{1i}$ is equal to $q_i(a_1,\ldots,a_n)$,
\[
v_{1i}=q_i(a_1,\ldots,a_n),\quad 1\le i\le r.
\]
\item The values of $q_i$ over $D$ are bounded by
the minimum and the maximum of
$\{v_{1i},\ldots,v_{si}\}$,
\[
\min(\{v_{ji}\}_{j=1}^s)\le q_i(\mathbf{x})\le \max(\{v_{ji}\}_{j=1}^s),\quad\forall\mathbf{x}\in D,\,1\le i\le r.
\]
\end{itemize}
\end{corollary}
\begin{proof}
Define the matrix $v=(v_1,\ldots,v_r)\in\mathbb{R}^{s\times r}$,
where $v_i=v_i(q_i)$ is defined as in Proposition \ref{mat-vec}.
\end{proof}

Combining the previous results on boxes and polynomials,
we obtain the following important result,

\begin{theorem}\label{main}
Let $\mathbf{d}=(d_1,\ldots,d_n)$ be a multi-degree,
$s=(d_1+1)\ldots(d_n+1)$ and
let $D=[a_1,a_1+2^{k_1}]\times\ldots\times[a_n,a_n+2^{k_n}]$
be a box, where $\mathbf{k}=(k_1,\ldots,k_n)\in\mathbb{N}_0^n$ is fixed.

For every $q_1,\ldots,q_r\in\mathbb{Z}[\mathbf{x}]_{\le\mathbf{d}}$
there exists a rectangular matrix $v\in\mathbb{Q}^{s\times r}$
and for every box $D'\subseteq D$ of the form
\[
D'=[z_1,z_1+2^{k_1'}]\times\ldots\times[z_n,z_n+2^{k_n'}],
\]
\[
z_i=a_i+l_i2^{k_i'},\,
0\le k'_i\le k_i,\,
0\le l_i<2^{k_i-k'_i},\,l_i\in\mathbb{N}_0,
0\le i\le n,
\]
there exists a
$(\mathbf{d},\mathbf{k})$-subdivision matrix $M\in\mathbb{Q}^{s\times s}$
such that $w=Mv$ satisfies
\begin{itemize}
\item The coefficient $w_{1i}$ is equal to $q_i(z_1,\ldots,z_n)$,
\[
w_{1i}=q_i(z_1,\ldots,z_n),\quad 1\le i\le r.
\]
\item The
values of $q_i$ over $D'$ are bounded by
the minimum and the maximum of
$\{w_{1i},\ldots,w_{si}\}$,
\[
\min(\{w_{ji}\}_{j=1}^s)\le q_i(\mathbf{x})\le \max(\{w_{ji}\}_{j=1}^s),\quad
1\le i\le r,\,
\mathbf{x}\in D'.
\]
\end{itemize}
\end{theorem}
\begin{proof}
Let $\varphi(x_1,\ldots,x_n)=(a_1+2^{k_1}x_1,\ldots,a_n+2^{k_n}x_n)$
and consider
the system $\varphi^{*}(q_1),\ldots,\varphi^{*}(q_r)$ over $[0,1]^n$.
Let $v\in\mathbb{Q}^{s\times r}$ be
the rectangular matrix such that its $i$-th column
is the vector representing $\varphi^{*}(q_i)$
in the Bernstein basis, $1\le i\le r$.

The idea now is to compare the system $\{q_i\}_{i=1}^r$
over $D$ and the system $\{\varphi^{*}(q_i)\}_{i=1}^r$ over $[0,1]^n$.
First of all, it is easy to check that $\varphi$ maps
bijectively the box $A\subseteq[0,1]^n$
to $D'\subseteq D$,
\[
A=\prod_{i=1}^n\left[\frac{z_i-a_i}{2^{k_i}},\frac{z_i-a_i+2^{k_i'}}{2^{k_i}}\right],
\quad
D'=\prod_{i=1}^n\left[z_i,z_i+2^{k_i'}\right].
\]
Let $\phi_A$ be the affine isomorphism sending $[0,1]^n$ to $A$
and let $\phi_{D'}$ the one sending $[0,1]^n$ to $D'$,
that is,
$\phi_A(x_1,\ldots,x_n)=(
(z_1-a_1+x_12^{k_1'})/2^{k_1},\ldots,
(z_n-a_n+x_n2^{k_n'})/2^{k_n})$
and
$\phi_{D'}(x_1,\ldots,x_n)=(
z_1+x_12^{k_1'},\ldots,z_n+x_n2^{k_n'})$.
Clearly $\phi_{D'}=\varphi\phi_A$, hence the
Bernstein coefficients of $\phi_{D'}^{*}(q_i)$
are the same as the
Bernstein coefficients of $\phi_A^{*}(\,\varphi^{*}(q_i)\,)$.
Then, the values of
$\{q_i\}_{i=1}^r$
over $D'$ are the same as the values of
$\{\varphi^{*}(q_i)\}_{i=1}^r$ over
$A$.
The result follows from Proposition \ref{subdiv} by
noting that
$A=\prod_{i=1}^n[l_i/2^{k_i-k_i'},(l_i+1)/2^{k_i-k_i'}]$
and letting $M$ be the matrix associated to
$A$ and $\mathbf{d}$.
\end{proof}

\begin{definition}\label{def-box}
Boxes appearing in the previous theorem
can be parameterized by two multi-indexes $(\mathbf{l},\mathbf{k'})$.
Let $\mathbf{a}\in\mathbb{Z}^n$, $\mathbf{k}\in\mathbb{N}_0^n$
and let $D=\prod_{i=1}^n[a_i,a_i+2^{k_i}]$.
For every $(\mathbf{l},\mathbf{k'})$
such that $\mathbf{0}\le\mathbf{k'}\le\mathbf{k}$ and
$0\le l_i<2^{k_i-k'_i}$, $l_i\in\mathbb{N}_0$, $1\le i\le n$
define
$D_{\mathbf{l},\mathbf{k'}}$ as
\[
D_{\mathbf{l},\mathbf{k'}}=\prod_{i=1}^n[a_i+l_i2^{k'_i},a_i+(l_i+1)2^{k'_i}].
\]
Otherwise, define $D_{\mathbf{l},\mathbf{k'}}=\emptyset$.
The vectors $\mathbf{a}$ and $\mathbf{k}$ are omitted from the notation.
Notice that $D_{\mathbf{0},\mathbf{k}}=D$
and that the volume of $D_{\mathbf{l},\mathbf{k'}}$
is $2^{k_1'+\ldots+k_n'}$ (or $0$). The index $\mathbf{l}$
represents the position and the index $\mathbf{k'}$
the volume of the box $D_{\mathbf{l},\mathbf{k'}}$.

\qed
\end{definition}

\section{Subdivision algorithm}\label{sect-2}

Let $\mathbf{d}\in\mathbb{N}_0^n$ be a multi-degree,
$s=(d_1+1)\ldots(d_n+1)$ and
let $\{M_i^L(\mathbf{d}),M_i^R(\mathbf{d})\}_{i=1}^n\subseteq\mathbb{Q}^{s\times s}$.
Given that these matrices depend only on
$\mathbf{d}$, it is possible to have them previously computed.
For example, depending on the application,
we can compute these matrices
to solve systems in $\mathbb{Z}[\mathbf{x}]_{\le\mathbf{d}}$
for some large $\mathbf{d}>>\mathbf{0}$.

Let $\mathbf{a}\in\mathbb{Z}^n$, $\mathbf{k}\in\mathbb{N}_0^n$,
$D=[a_1,a_1+2^{k_1}]\times\ldots\times[a_n,a_n+2^{k_n}]$
and let $q_1,\ldots,q_r\in\mathbb{Z}[\mathbf{x}]_{\le\mathbf{d}}$.
We want to describe the following sets,
\[
\mathcal{S'}=\left\{\mathbf{z}\in D\cap\mathbb{Z}^n\,\colon\,
q_2(\mathbf{z})\ge 0,\,\ldots\,,q_r(\mathbf{z})\ge 0\right\},
\quad
\mathcal{S}=\left\{\mathbf{z}\in \mathcal{S'}\,\colon
q_1(\mathbf{z})=\max_{\mathcal{S'}}(q_1)\right\}.
\]
Let $v\in\mathbb{Q}^{s\times r}$ denote
the rectangular matrix associated to $q_1,\ldots,q_r$ (see Theorem \ref{main}).
Initialize the global variables $\mathcal{S}:=\emptyset$
and $\theta:=-\infty$
and run the function
{\tt SolveIPP}($v$, $\mathbf{0}$, $\mathbf{k}$).

\begin{function}[H]
\DontPrintSemicolon
\caption{SolveIPP($w$, $\mathbf{l}$, $\mathbf{k'}$)}
\KwData{Global variable $\theta$, a global set $\mathcal{S}$
 and matrices $\{M_i^L(\mathbf{d}),M_i^R(\mathbf{d})\}_{i=1}^n$.}
\BlankLine
\If{$\max(w_1) \ge\theta$ {\bf and} $\max(w_2)\ge 0$ {\bf and} $\ldots$ {\bf and} $\max(w_r)\ge 0$}{\label{ste1}
 $k_j':=\max(\mathbf{k'})$\;
 \If{$k_j'>0$}{\label{ste3}
  \SolveIPP($M_j^Lw$, $\mathbf{l}+l_je_j$, $\mathbf{k'}-e_j$)\;
  \SolveIPP($M_j^Rw$, $\mathbf{l}+l_je_j+e_j$, $\mathbf{k'}-e_j$)\;
 }
 \ElseIf{$w_{11}\ge\theta$ {\bf and} $w_{12}\ge 0$ {\bf and}
   $\ldots$ {\bf and} $w_{1r}\ge 0$}{\label{ste6}
  \If{$w_{11}>\theta$}{
   $\mathcal{S}:=\{\}$\;
   $\theta:=w_{11}$\;
  }
  $\mathcal{S}:=\mathcal{S}\cup\{\mathbf{a}+\mathbf{l}\}$\;
 }
}
\end{function}

Let us explain the subdivision process of \SolveIPP (Solve Integer Polynomial Programming). 
The first thing to say
is that it is a recursive algorithm that makes
two matrix multiplications in each call, $M_j^Lv$ and $M_j^Rv$.
The rationality of the matrices makes this algorithm robust and reliable
and the resulting set $\mathcal{S}$
is equal to its mathematical counterpart $\mathcal{S}$ (see Theorem \ref{answer} below).

According to Theorem \ref{main} the values in
$(w,\mathbf{l},\mathbf{k'})$ give information
about the polynomials $\{q_1,\ldots,q_r\}$
over the box $D_{\mathbf{l},\mathbf{k'}}$.
The condition in Step \ref{ste1} is given to
test if there exist solutions on $D_{\mathbf{l},\mathbf{k'}}$,
in other words, if we keep $D_{\mathbf{l},\mathbf{k'}}$
or we reject it.
Assuming that we do not reject $D_{\mathbf{l},\mathbf{k'}}$,
we test in Step \ref{ste3} if
$D_{\mathbf{l},\mathbf{k'}}$ can be divided.

If $D_{\mathbf{l},\mathbf{k'}}$ can be divided,
we define $D_L=D_{\mathbf{l}+l_je_j,\mathbf{k'}-e_j}$
and $D_R=D_{\mathbf{l}+l_je_j+e_j,\mathbf{k'}-e_j}$
and we analyze the values of the system over
these boxes by defining $w^L=M_j^Lv$ and $w^R=M_j^Rv$.
The notation $e_j$ indicates the $j$-th cannonical vector,
\[
(e_j)_i=
\begin{cases}
1& i=j\\
0&i\neq j
\end{cases},\quad 1\le i\le n.
\]

If $D_{\mathbf{l},\mathbf{k'}}$ can not be divided,
that is $\mathbf{k'}=\mathbf{0}$,
we test in Step \ref{ste6} if $\mathbf{a}+\mathbf{l}$
is actually a solution of the system.
If this is the case, we save $\mathbf{z}=\mathbf{a}+\mathbf{l}$
in $\mathcal{S}$.
If $w_{11}>\theta$, we reset
$\mathcal{S}$ and set $\theta=q_1(\mathbf{z})$.
This last process guarantees that
$q_1(\mathbf{z})=\theta$ for every $\mathbf{z}\in \mathcal{S}$.

\begin{theorem}\label{answer}
Let $n\in\mathbb{N}$, $\mathbf{d},\mathbf{k}\in\mathbb{N}_0^n$,
$\mathbf{a}\in\mathbb{Z}^n$,
$q_1,\ldots,q_r\in\mathbb{Z}[x_1,\ldots,x_n]_{\le\mathbf{d}}$,
let $D=\prod_{i=1}^n[a_i,a_i+2^{k_i}]$ be a box
and
let $v$ be the matrix associated to $q_1,\ldots,q_r$ as in Theorem \ref{main}.
Then, the output of the 
function \SolveIPP($v$, $\mathbf{0}$, $\mathbf{k}$) is the set $\mathcal{S}$,
\[
\mathcal{S'}=\left\{\mathbf{z}\in D\cap\mathbb{Z}^n\,\colon\,
q_2(\mathbf{z})\ge 0,\,\ldots\,,q_r(\mathbf{z})\ge 0\right\},
\quad
\mathcal{S}=\left\{\mathbf{z}\in \mathcal{S'}\,\colon
q_1(\mathbf{z})=\max_{\mathcal{S'}}(q_1)\right\}.
\]
\end{theorem}
\begin{proof}
The compactness of $D$ implies that the set $\mathcal{S'}$ is finite.
If $\mathcal{S'}\neq\emptyset$,
the set $q_1(\mathcal{S'})$ has a maximum value $\theta$
and the space of solutions $\mathcal{S}$ is equal to
$q_1^{-1}(\theta)\cap\mathcal{S'}$.

Let us prove that every $\mathbf{z}\in\mathcal{S}$
is in the output of \SolveIPP (the reciprocal is straightforward).
Let $D'$ be a box
such that $\mathbf{z}\in D'\cap \mathcal{S}$. Then
$\max_{D'}(q_1)\ge \theta$ and
$\max_{D'}(q_i)\ge 0$, $2\le i\le r$.
In particular,
the box $D'$ is not rejected by the algorithm
and is subdivided. Then, without
loss of generality, we may suppose that
$D'=\prod_{i=1}^n[z_i,z_i+1]$ and clearly,
$\mathbf{z}$ is in the output of \SolveIPP.
Note that after some steps,
the value of the global variable $\theta$
is equal to the maximum of $q_1$ over $\mathcal{S'}$.
\end{proof}

\section{Expected complexity}
Given that the number of matrix multiplications
can be computed as the number of times
the function \SolveIPP is called (minus one),
the complexity of the algorithm can be computed
in terms of boxes.
The worst case scenario is
when every box is subdivided (exponential case).
The best case occurs when almost every box is unnecessary (linear case).
Thus, the complexity of the algorithm
is related to the number of unnecessary boxes, and it can be described 
in probabilistic terms. The following proposition gives a geometric 
representation of the set of unnecessary boxes.

\begin{proposition}\label{elev2}
Let $n\in\mathbb{N}$, $\mathbf{d},\mathbf{k}\in\mathbb{N}_0^n$,
$\mathbf{a}\in\mathbb{Z}^n$,
$q_1,\ldots,q_r\in\mathbb{Z}[x_1,\ldots,x_n]_{\le\mathbf{d}}$,
$D=\prod_{i=1}^n[a_i,a_i+2^{k_i}]$ and
let $\theta$ be the maximum of $q_1$ over
$\mathcal{S'}=\{q_i\ge 0\}_{i=2}^r\cap D\cap\mathbb{Z}^n$.
Let $\mathcal{R}$ be the set
of real solutions of the system,
\[
\mathcal{R}=\{\mathbf{z}\in D
\,\colon\, q_1(\mathbf{z})\ge\theta,\,q_{2}(\mathbf{z})\ge 0,\ldots,
q_r(\mathbf{z})\ge 0\}.
\]
Then, there exists
$\mathbf{d'}\in\mathbb{N}_0^n$
such that for every $\mathbf{d''}\ge\mathbf{d'}$
the following sentences are equivalent for every
non-empty box $D_{\mathbf{l},\mathbf{k'}}$:
\begin{enumerate}
\item The matrix associated to the system
over $D_{\mathbf{l},\mathbf{k'}}$ satisfies the condition in Step \ref{ste1}.
\item $D_{\mathbf{l},\mathbf{k'}}\cap\mathcal{R}\neq\emptyset$.
\item There exists $D_{\mathbf{l'},\mathbf{0}}\subseteq
D_{\mathbf{l},\mathbf{k'}}$
such that $D_{\mathbf{l'},\mathbf{0}}\cap\mathcal{R}\neq\emptyset$.
\end{enumerate}
\end{proposition}
\begin{proof}
Let us prove the existence of $\mathbf{d'}$.
First, replace $q_1$ by $q_1-\theta$, and write
condition in Step \ref{ste1} as $\{\max(w_i)\ge 0\}_{i=1}^r$.
In each box $D_{\mathbf{l},\mathbf{k'}}$
take some $\mathbf{d}_{\mathbf{l},\mathbf{k'}}$ satisfying
\[
\max_{D_{\mathbf{l},\mathbf{k'}}}(q_i)<0\Longrightarrow\max(w_i)<0,\quad 1\le i\le r,
\]
where the rectangular matrix $w=(w_1,w_2,\ldots,w_r)$
is the one associated to the system
over the box $D_{\mathbf{l},\mathbf{k'}}$
and multi-degree $\mathbf{d}_{\mathbf{l},\mathbf{k'}}$,
see Lemma \ref{elev}.
Let $\mathbf{d'}$ be the maximum of
$\{\mathbf{d}_{\mathbf{l},\mathbf{k'}}\}$.

Assume that we have a box $D_{\mathbf{l},\mathbf{k'}}$
without real solutions in it. Let
us see that the algorithm applied in multi-degree $\mathbf{d'}$
rejects this box immediately.
Given that $D_{\mathbf{l},\mathbf{k'}}$ does not
contain real solutions, there
exists some $i$, $1\le i\le r$, such that $\min(q_i)<0$.
Then, $\min(w_{i})<0$ and the box $D_{\mathbf{l},\mathbf{k'}}$ is rejected.

The equivalence between the last two sentences follows
from the fact that it is possible to cover
$D_{\mathbf{l},\mathbf{k'}}$ with boxes of the
form $D_{\mathbf{l'},\mathbf{0}}$.
\end{proof}

\begin{definition}\label{complex-number}
Let $n\in\mathbb{N}$, $\mathbf{d},\mathbf{k}\in\mathbb{N}_0^n$,
$\mathbf{a}\in\mathbb{Z}^n$,
$q_1,\ldots,q_r\in\mathbb{Z}[x_1,\ldots,x_n]_{\le\mathbf{d}}$
and $D=\prod_{i=1}^n[a_i,a_i+2^{k_i}]$,
let $\theta$ be the maximum of $q_1$ over
$\mathcal{S'}=\{q_i\ge 0\}_{i=2}^r\cap D\cap\mathbb{Z}^n$
and let $\mathcal{R}$ be the set
of real solutions,
$\mathcal{R}=\{\mathbf{z}\in D
\,\colon\, q_1(\mathbf{z})\ge\theta,\,q_{2}(\mathbf{z})\ge 0,\ldots,
q_r(\mathbf{z})\ge 0\}$.

The \emph{complexity number} of the system $\{q_1,\ldots,q_r\}$
over $D$ is defined as
\[
\lambda(\mathcal{R})=\frac{\#\{\mathbf{l}\,\colon\,
D_{\mathbf{l},\mathbf{0}}\cap\mathcal{R}\neq\emptyset\}}{2^{k_1+\ldots+k_n}},
\quad \lambda(\mathcal{R})\in\mathbb{Q},\,0\le \lambda(\mathcal{R})\le 1.
\]
The complexity number can be described in probabilistic terms,
it is the probability of a
box $D_{\mathbf{l},\mathbf{0}}$ to have real solutions.
\qed
\end{definition}

\begin{theorem}\label{comple}
Let $\lambda$ be the complexity number of
$\{q_1,\ldots,q_r\}$ over $D=\prod_{i=1}^n[a_i,a_i+2^{k_i}]$
and let $K=k_1+\ldots+k_n$.
Let $\mathbf{d'}$ be
the multi-degree from Proposition \ref{elev2}.
Then, the expected number
of multiplications between a square triangular
matrix in $\mathbb{Q}^{s'\times s'}$
and a rectangular matrix in $\mathbb{Q}^{s'\times r}$,
$s'=(d_1'+1)\ldots(d_n'+1)$ produced by \SolveIPP is
\[
\sum_{i=0}^K
2^{i}\prod_{j=0}^{i-1}\left(1-(1-\lambda)^{2^{K-j}}\right)-1.
\]
This number is bounded between
$((2\lambda)^{K+1}-1)/(2\lambda-1)-1$ (or $K$ if $\lambda=1/2$)
and $2(2^{K}-1)$.
\end{theorem}
\begin{proof}
Each time a box is subdivided, the function \SolveIPP
produces two matrix multiplications.
Hence, the expected number of matrix multiplications
is equal to the expected number of divided boxes.
Equivalently, it is the expected number of boxes
processed by \SolveIPP minus 1 ($D$ is not
a divided box).
According to Proposition \ref{elev2}, the probability
of a box $D_{\mathbf{l},\mathbf{k'}}$ to be subdivided
depends on the existence of some $D_{\mathbf{l'},\mathbf{0}}$
such that $D_{\mathbf{l'},\mathbf{0}}
\subseteq D_{\mathbf{l},\mathbf{k'}}$ and $D_{\mathbf{l'},\mathbf{0}}\cap\mathcal{R}\neq 0$.
There are a total of $2^{k_1'+\ldots+k_n'}$ possible boxes
$D_{\mathbf{l'},\mathbf{0}}$ inside $D_{\mathbf{l},\mathbf{k'}}$
hence the probability
is equal to $1-(1-\lambda)^{2^{k_1'+\ldots+k_n'}}$.
We can model the problem of finding the expected number of boxes
as a \emph{Branching process}, see \citep{MR0163361}.
We say that $D_{\mathbf{l},\mathbf{k'}}$
is in \emph{generation} $i$
if $k_1'+\ldots+k_n'=K-i$.
In other words, if the volume of $D_{\mathbf{l},\mathbf{k'}}$
is equal to $2^{K-i}$.
Let $\{Z_i\}_{i=0}^{K}$ be the
number of boxes in each generation.
Boxes in generation $i$
have probability $\lambda_i$ (resp. $1-\lambda_i$)
to add two (resp. $0$) boxes to the next generation,
\[
\lambda_i=1-(1-\lambda)^{2^{K-i}},\quad 0\le i\le K.
\]
These boxes can be interpreted as
independent identically distributed random variables $\{D_i\}$
with common generating function $1-\lambda_i+\lambda_ix^2$.

Thus $Z_{i+1}=D_1+\ldots+D_{Z_i}$ and
from \cite[\S 5.1 Th.25]{MR2059709} we obtain
$E(Z_{i+1})=2\lambda_iE(Z_{i})$, hence
\[
E(Z_0)=1,\quad
E(Z_{i+1})=
2\lambda_iE(Z_{i})=2^{i+1}\lambda_i\ldots\lambda_0,\quad 0\le i\le K-1.
\]
The number of expected boxes is
$E(Z_0)+\ldots+E(Z_K)$.
\end{proof}

\begin{example*}
Consider the equation $y=x$ in $[0,2]\times[0,2]$.
Clearly, it
has two boxes with real solutions,
$[0,1]\times[0,1]$ and $[1,2]\times[1,2]$. Then, $\lambda=2/4$
and the expected number of boxes processed by \SolveIPP
for $\mathbf{d'}$ as in Proposition \ref{elev2}
is $91/16$, a number
between $2+1=3$ and $2^3-1=7$.
The function applied to this example
passes Step \ref{ste1}, $7$ times.
\qed
\end{example*}

\begin{example*}
Consider the system $y-x^2=0$ in $D=[0,2^3]\times [0,2^3]$.
The system has three integral solutions, $(0,0)$, $(1,1)$
and $(2,4)$. The real curve pass through $8$ boxes
of the form $D_{\mathbf{l},\mathbf{0}}$.
Then,
\[
\lambda=8/64,\quad (\lambda_0,\ldots,\lambda_7)\cong
(.99, .98, .88, .65, .41, .23, .12).
\]
Hence, the expected number of boxes processed by \SolveIPP
for $\mathbf{d'}$ as in Proposition \ref{elev2}
is approximately $34$.
In this example Step \ref{ste1} occurs $33$ times.
\qed
\end{example*}

\begin{remark*}
In \citep{patil2014improved} the authors manage to 
to improve (a variant of) \SolveIPP
to produce better performance. Also, it is 
possible to adapt \SolveIPP to treat
systems of continuous functions
by using Stone-Weierstrass
Theorem \citep[Cor.1]{MR0027121} and interpolation
techniques \citep{MR936450,MR864015}.

Notice that Proposition \ref{elev2} implies that the strategy of
dividing a side of a box in half is not optimal.
For example, consider the equation $x-2$ over $[0,4]$.
The function \SolveIPP divides $[0,4]$ in two halves
making two matrix multiplications,
but given that \SolveIPP do not reject these two halves, 
it makes four multiplications more.

Instead of dividing $[0,4]$ into $[0,2]$ and $[2,4]$
it is better (from a performance point of view) 
to divide it into $[0,2)$ and $[2,4]$
or even better into $[0,1]$ and $[2,4]$.
\qed
\end{remark*}

\begin{remark*}
Given $\mathbf{d}\in\mathbb{N}^n$,
we need to compute the triangular matrices
$\{M_i^L,M_i^R\}_{i=1}^n\subseteq\mathbb{Q}^{s\times s}$,
$s=\prod_{i=1}^n(d_i+1)$ and save them into the computer memory.
If $ns(s+1)$ is greater than the total amount of available memory
it is possible to use techniques from computer science such as
Parallel computing or Map-Reduce
\citep{Rajaraman:2011:MMD:2124405,Dean:2008:MSD:1327452.1327492}
in order to increase the available space of memory.
An alternative (or complementary) 
solution is to use \emph{toric} Bernstein polynomials to minimize
the number of monomials, \citep[Cor.2]{MR2496414}.
\qed
\end{remark*}

\section{Solving systems of congruence equations}

In this section we apply the subdivision method presented in Section \ref{sect-2}
to solve congruence equations ({\tt SolveCE}).
Let $h_1,\ldots,h_r\in \mathbb{Z}[x_1,\ldots,x_n]_{\le\mathbf{d}}$
and let $p_1,\ldots,p_r\in\mathbb{Z}_{\ge 2}$ be natural numbers (possibly coprimes). Here we propose
a method to find solutions $\mathbf{z}\in \mathbb{Z}^{n}$ to the system
\[
\left\{
\begin{array}{lcrr}
h_1(\mathbf{z})&\equiv&0, & \mbox{mod}\,(p_1),\\
&\vdots\\
h_r(\mathbf{z})&\equiv&0, & \mbox{mod}\,(p_r).
\end{array}
\right.
\]
In order to do so, we need to adapt some of the propositions presented so far.
Notice that any solution $\mathbf{z}\in \mathbb{Z}^{n}$
can be represented in the box $[0,2^k]\times\ldots\times[0,2^k]$, where
\[
k:=\lceil\log_2(\mbox{lcm}(p_1,\ldots,p_r))\rceil \in\mathbb{Z}_{\ge 2},
\]
$\mbox{lcm}(x,y)$ is the least common multiple of $x$ and $y$
and $\lceil x\rceil$ is the least integer greater than or equal to $x$.
The following simple Lemma is the key ingredient of our method.
\begin{lemma}\label{lemma-hasSol}
Let $p\in\mathbb{Z}_{\ge 2}$, let $m_1\le m_2$ be
two integers and let $r_1$ (resp. $r_2$) be
the remainder of $m_1$ (resp. $m_2$) divided by $p$.
The following sentences are equivalent,
\begin{enumerate}
\item\label{i1} $m_2-m_1<p-1$ and $0<r_1\le r_2$.
\item\label{i2} There is no multiple of $p$ between $m_1$ and $m_2$.
\end{enumerate}
\end{lemma}
\begin{proof}
Let $q_1,q_2\in\mathbb{Z}$ be such that $m_1=q_1p+r_1$
and $m_2=q_2p+r_2$. Given that $m_1\le m_2<m_1+p-1$, we have
the inequality,
\[
p-1>m_2-m_1=(q_2-q_1)p+(r_2-r_1)\ge0.
\]

\noindent \ref{i1}$\Rightarrow$\ref{i2})
Given that $0\le r_2-r_1\le p-1$, we have $(q_2-q_1)p\ge 0$
and from $p-1>(q_2-q_1)p\ge0$ it follows $q_2=q_1$.
Then, there exists $q\in\mathbb{Z}$ such that
$m_1=qp+r_1$ and $m_2=qp+r_2$.
Assume now, that there exists $a\in\mathbb{Z}$ such that
$m_1\le ap\le m_2$.
Then, $m_2-ap$ is also $<p-1$ and
it follows
$a=q$. Hence,
\[
m_1\le qp\le m_2\iff
qp+r_1\le qp\le qp+r_2\iff r_1\le0\le r_2.
\]
A contradiction, because $r_1>0$.

\

\noindent \ref{i2}$\Rightarrow$\ref{i1})
It is clear that if
$m_2-m_1\ge p-1$, then there exists a multiple
of $p$ between $m_1$ and $m_2$. Same if $r_1=0$ or $r_2=0$.
Assume $m_2-m_1<p-1$ and $r_1>r_2$
and let us prove that there exists a multiple
of $p$ between $m_1$ and $m_2$.

Using the inequality
$p-1>(q_2-q_1)p+(r_2-r_1)\ge0$ and $0<r_1-r_2\le p-1$
it follows that $q_2=q_1+1$. Then, $q_1p$ is in between
$m_1$ and $m_2$.
\end{proof}

A direct application of Lemma \ref{lemma-hasSol}
gives the following function,

\begin{function}[H]
\DontPrintSemicolon
\caption{hasSol($m_1$,$m_2$,$p$)}
\KwData{$m_1\le m_2\in\mathbb{Q}$ and $p\in\mathbb{Z}_{\ge 2}$.}
\KwOut{True if there is a multiple of $p$ between $m_1$ and $m_2$. False otherwise.}
\BlankLine
$(m_1,m_2)\leftarrow (\lceil m_1\rceil,\lfloor m_2\rfloor)$\;
$(r_1,r_2)\leftarrow (m_1 \% p,m_2 \% p)$\;
\If{$m_2-m_1<p-1$ {\bf and} $0<r_1\le r_2$}{
\KwRet{\tt False}
}
\KwRet{\tt True}
\end{function}

\noindent Step 1 in \hasSol is based on
the following fact. There exists an
integer $i$ such that $m_1\le i\le m_2$
if only if $\lceil m_1\rceil\le i\le \lfloor m_2\rfloor$,
where $\lceil m_1\rceil$
is the least integer greater than or equal to $m_1$
and $\lfloor m_2\rfloor$
is the greatest integer less than or equal to $m_2$.
The notation $m_1\% p$ in Step 2 means the \emph{remainder}
of dividing $m_1$ by $p$. Same for $m_2\% p$.

\

Now, we can adapt the function \SolveIPP.
Fix some $\mathbf{d}\in\mathbb{N}_0^n$ and
let $\{M_i^L(\mathbf{d}),M_i^R(\mathbf{d})\}_{i=1}^n$
be the matrices defined in \ref{def-m}.
Let $h_1,\ldots,h_r\in \mathbb{Z}[\mathbf{x}]_{\le\mathbf{d}}$
and let $p_1,\ldots,p_r\in\mathbb{Z}_{\ge 2}$ be natural numbers.
Let $k:=\lceil\log_2(\mbox{lcm}(p_1,\ldots,p_r))\rceil \in\mathbb{Z}_{\ge 2}$
and $D=[0,2^{k}]\times\ldots\times[0,2^{k}]$. Denote $\mathbf{p}=(p_1,\ldots,p_r)$,
$\mathbf{0}=(0,\ldots,0)$ and $\mathbf{k}=(k,\ldots,k)$.

We want to describe the following set,
\[
\mathcal{S'}=\left\{\mathbf{z}\in D\cap \mathbb{Z}^{n}\,\colon\,
h_1(\mathbf{z})\equiv 0\mbox{ mod }(p_1),\ldots,
h_r(\mathbf{z})\equiv 0\mbox{ mod }(p_r)\right\}.
\]
Let $v$ be the rectangular matrix associated to $h_1,\ldots,h_r$
as in Theorem \ref{main}.
Initialize the global set $\mathcal{S'}:=\emptyset$
and run {\tt SolveCE}($\mathbf{p}$,$v$, $\mathbf{k}$,$\mathbf{0}$).

\begin{function}[H]
\DontPrintSemicolon
\caption{SolveCE($\mathbf{p}$,$w$,$\mathbf{k'}$, $\mathbf{l}$)}
\KwData{A global set $\mathcal{S'}$
 and matrices $\{M_i^L(\mathbf{d}),M_i^R(\mathbf{d})\}_{i=1}^n$.}
\KwOut{$\mathcal{S'}=
\left\{\mathbf{z}\in D\cap\mathbb{Z}^{n}\,\colon\,
h_1(\mathbf{z})\equiv 0\mbox{ mod }(p_1),\ldots,
h_r(\mathbf{z})\equiv 0\mbox{ mod }(p_r)\right\}
$}
\BlankLine
\If{\hasSol($\min(w_1),\max(w_1),p_1$) {\bf and} $\ldots$ {\bf and}
\hasSol($\min(w_r),\max(w_r),p_r)$}{
 $k_j':=\max(\mathbf{k'})$\;
 \If{$k_j'>0$}{
  \SolveCE($\mathbf{p}$, $M_j^Lw$, $\mathbf{k'}-e_j$, $\mathbf{l}+l_je_j$)\;
  \SolveCE($\mathbf{p}$, $M_j^Rw$, $\mathbf{k'}-e_j$, $\mathbf{l}+l_je_j+e_j$)\;
 }
 \ElseIf{$p_1|w_{11}$ {\bf and} $\ldots$ {\bf and} $p_r|w_{1r}$}{
  $\mathcal{S'}:=\mathcal{S'}\cup\{\mathbf{l}\}$\;
 }
}
\end{function}

\begin{theorem}\label{answer2}
The output of \SolveCE($\mathbf{p}$,$v$, $(k,\ldots,k)$,$(0,\ldots,0)$)
is the set
\[
\mathcal{S'}=
\left\{\mathbf{z}\in D\cap\mathbb{Z}^{n}\,\colon\,
p_1|h_1(\mathbf{z}),\ldots,
p_r|h_r(\mathbf{z})\right\},
\]
where
$D=[0,2^{k}]\times\ldots\times[0,2^{k}]$,
$k=\lceil\log_2(\mbox{lcm}(p_1,\ldots,p_r))\rceil$ and
$v$ is the matrix associated to $h_1,\ldots,h_r$
as in Theorem \ref{main}.
\end{theorem}
\begin{proof}
If there is a solution in a box $D'$,
then the process goes from Step 1 to Step 2
and $D'$ gets subdivided unless the volume of $D'$ is equal to $1$.
If the volume of $D'$ is equal to $1$ and there exists
a solution $\mathbf{l}\in D'$,
then the condition in Step 6 is satisfied and
$\mathbf{l}$ is saved in $\mathcal{S'}$.
\end{proof}

As before, we can estimate the expected number of
matrix multiplications that \SolveCE do. Let
$\mathbf{h}:\mathbb{R}^n\to\mathbb{R}^r$ be the polynomial
function $\mathbf{h}=(h_1,\ldots,h_r)$ and consider the
discrete set
\[
\mathcal{X}=\{(x_1,\ldots,x_r)\in\mathbb{Z}^r\,\colon\, p_1|x_1,\ldots,p_r|x_r\}.
\]
Notice that the set $\mathbf{h}^{-1}(x)\subseteq\mathbb{R}^n$ is a real
algebraic variety for each $x\in\mathcal{X}$.

\begin{proposition}\label{elev3}
Let $n\in\mathbb{N}$, $\mathbf{d}\in\mathbb{N}_0^n$,
$h_1,\ldots,h_r\in\mathbb{Z}[x_1,\ldots,x_n]_{\le\mathbf{d}}$,
$p_1,\ldots,p_r\in\mathbb{Z}_{\ge 2}$
and let
$k=\lceil\log_2(\mbox{lcm}(p_1,\ldots,p_r))\rceil$.
Let
$\mathcal{X}=\{\mathbf{x}\in\mathbb{Z}^r\,\colon\, p_1|x_1,\ldots,p_r|x_r\}$
and let $\mathbf{h}:\mathbb{R}^n\to\mathbb{R}^r$, $\mathbf{h}=(h_1,\ldots,h_r)$.
Then, there exists
$\mathbf{d'}\in\mathbb{N}_0^n$
such that for every $\mathbf{d''}\ge\mathbf{d'}$
and every non-empty box $D_{\mathbf{l},\mathbf{k'}}$,
the following sentences are equivalent:
\begin{enumerate}
\item The matrix associated to $(h_1,\ldots,h_r)$
over $D_{\mathbf{l},\mathbf{k'}}$ satisfies the condition in Step 1.
\item $D_{\mathbf{l},\mathbf{k'}}\cap\mathbf{h}^{-1}(\mathcal{X})\neq\emptyset$.
\item There exists $D_{\mathbf{l'},\mathbf{0}}\subseteq
D_{\mathbf{l},\mathbf{k'}}$
such that $D_{\mathbf{l'},\mathbf{0}}\cap\mathbf{h}^{-1}(\mathcal{X})\neq\emptyset$.
\end{enumerate}
\end{proposition}
\begin{proof}
Clearly sentence 2 is equivalent to sentence 3. Let us prove
the equivalence between sentences 1 and 2.
Fix some $\mathbf{l}$ and $\mathbf{k'}$.
Let us define $m_1^i$ and $m_2^i$, $1\le i\le r$,
\[
m_1^i:=
\begin{cases}
\lfloor\min_{D_{\mathbf{l},\mathbf{k'}}}(h_i)\rfloor
&\mbox{if }
\min_{D_{\mathbf{l},\mathbf{k'}}}(h_i)
>\lfloor\min_{D_{\mathbf{l},\mathbf{k'}}}(h_i)\rfloor\\
\min_{D_{\mathbf{l},\mathbf{k'}}}(h_i)-\frac{1}{2}
&\mbox{if not}
\end{cases}
\]
\[
m_2^i:=
\begin{cases}
\lceil\max_{D_{\mathbf{l},\mathbf{k'}}}(h_i)\rceil
&\mbox{if }
\max_{D_{\mathbf{l},\mathbf{k'}}}(h_i)
<\lceil\max_{D_{\mathbf{l},\mathbf{k'}}}(h_i)\rceil\\
\max_{D_{\mathbf{l},\mathbf{k'}}}(h_i)+\frac{1}{2}
&\mbox{if not}
\end{cases}
\]
Choose a multi-degree $\mathbf{d}_{\mathbf{l},\mathbf{k'}}$ as in Lemma \ref{elev}
such that
\[
m_1^i<\min(w_i)\le \max(w_i)<m_2^i,\quad 1\le i\le r,
\]
where $w=(w_1,\ldots,w_r)$ is the rectangular matrix
associated to $(h_1,\ldots,h_r)$ over $D_{\mathbf{l},\mathbf{k'}}$.
By construction, there exists a multiple of $p_i$ between
$\min_{D_{\mathbf{l},\mathbf{k'}}}(h_i)$ and
$\max_{D_{\mathbf{l},\mathbf{k'}}}(h_i)$ if and only if
it is between
$\min(w_i)$ and $\max(w_i)$. The result follows by taking
$\mathbf{d'}$ as the maximum of
$\{\mathbf{d}_{\mathbf{l},\mathbf{k'}}\}$.
\end{proof}

\begin{theorem}\label{comple2}
Let $h_1,\ldots,h_r\in\mathbb{Z}[x_1,\ldots,x_n]$,
$p_1,\ldots,p_r\in\mathbb{Z}_{\ge 2}$, 
$k=\lceil\log_2(\mbox{lcm}(p_1,\ldots,p_r))\rceil$
and let $\mathbf{d'}$ be as in Proposition \ref{elev3}.
Let $\mathcal{X}=\{\mathbf{x}\in\mathbb{Z}^r\,\colon\, p_1|x_1,\ldots,p_r|x_r\}$
and let $\mathbf{h}:\mathbb{R}^n\to\mathbb{R}^r$, $\mathbf{h}=(h_1,\ldots,h_r)$.
Then, the expected number of matrix multiplications produced by {\tt SolveCE}
in $D=\prod_{i=1}^n[0,2^{k}]$ is
\[
\sum_{i=0}^{nk}
2^{i}\prod_{j=0}^{i-1}\left(1-(1-\lambda)^{2^{nk-j}}\right)-1,\quad
\lambda:=\frac{\#\{\mathbf{l}\,\colon\,
D_{\mathbf{l},\mathbf{0}}\cap\mathbf{h}^{-1}(\mathcal{X})\neq\emptyset
\}}{2^{nk}}.
\]
The number $\lambda$ is the probability of a box $D_{\mathbf{l},\mathbf{0}}$ to intersect 
$\mathbf{h}^{-1}(\mathcal{X})$.
\end{theorem}
\begin{proof}
According to Proposition \ref{elev3}
the probability of a box $D_{\mathbf{l},\mathbf{k'}}$ to be rejected
in Step 1 is $(1-\lambda)^{2^{k_1'+\ldots+k_n'}}$.
The proof continues as in Theorem \ref{comple}.
\end{proof}

\begin{definition}\label{bf}
We define the \emph{brute-force algorithm} as
the algorithm consisting in testing
each $\mathbf{z}\in D\cap\mathbb{Z}^{n}$ to be in $\mathcal{S'}$,
\[
\mathcal{S'}=
\left\{\mathbf{z}\in D\cap\mathbb{Z}^{n}\,\colon\,
p_1|h_1(\mathbf{z}),\ldots,
p_r|h_r(\mathbf{z})\right\}.
\]
This method performs $r2^{nk}$ evaluations
to characterize $\mathcal{S'}$. Notice that there is
no need for the sides of $D$ to be a power of two
and the same is true in \SolveCE.
We chose to have sides to be a power of two
to make the presentation clearer and also to be
able to compare this algorithm with \SolveCE.

For example,
if $p_1=p_2=5$ and $h_1=x+1$, $h_2=y$ we need to take
the starting box to be $D=[0,8]^2$. The brute-force algorithm
with the box $D=[0,5]^2$ performs $25\times 2$ evaluations to
find the solution $(4,0)$.
The function \SolveCE finds two (apparently different)
solutions $(4,0)$ and $(4,5)$ in the box $D=[0,8]^2$.
\qed
\end{definition}

In order to avoid benefiting either algorithm, we choose
powers of two for the numbers $p_1,\ldots,p_r$.
In the next table we compare \SolveCE with
the brute-force algorithm. All computations were made
using our implementation in SageMath.
\begin{center}
\begin{tabular}{|c|c|c|c|c|}
\hline
$\mathbf{p}$&$\mathbf{h}$&Step 1&Step 6&BF\\
\hline\hline
$512$&$2x-3$&$19$&$2$&$512$\\
\hline
$64$&$x^2+3x-4$&$111$&$51$&$64$\\
\hline
$128$&$x^7+5x^2-9$&$254$&$127$&$128$\\
\hline
$(512,512)$&$(x+1,y)$&$40$&$4$&$262144$\\
\hline
$(8,8)$&$(xy+2x+y^2,3y+x^2)$&$123$&$60$&$64$\\
\hline
$(128,128)$&$(y(x-2),y^5+1)$&$29124$&$13675$&$16384$\\
\hline
$(16,16,16)$&$(x+y-z,2x+y-3z,z-3+x)$&$491$&$84$&$4096$\\
\hline
$(4,8,16)$&$(x+5yz-z^3,y^2,x+z+2)$&$1923$&$720$&$512$\\
\hline
\end{tabular}
\end{center}
The first two columns codify the system $\mathbf{h}\equiv \mathbf{0}$ mod $(\mathbf{p})$.
The column \emph{Step 1} is the number of times \SolveCE
pass through the first step.
The column \emph{Step 6} is the number of boxes $D_{\mathbf{l},\mathbf{0}}$
tested to have a solution. Finally, the last column is the number
of evaluations that the brute-force algorithm must do.
We can infer from the table that \SolveCE is faster for
linear systems, but for linear systems there are much
better algorithms, \citep{MR1689189,MR2312337}.
We have prioritized presentation over performance.


\begin{thebibliography}{10}

\bibitem{MR1992831}
Alexander Barvinok and Kevin Woods.
\newblock Short rational generating functions for lattice point problems.
\newblock {\em J. Amer. Math. Soc.}, 16(4):957--979 (electronic), 2003.

\bibitem{MR3473492}
Fani Boukouvala, Ruth Misener, and Christodoulos~A. Floudas.
\newblock Global optimization advances in {M}ixed-{I}nteger {N}onlinear
  {P}rogramming, {MINLP}, and {C}onstrained {D}erivative-{F}ree {O}ptimization,
  {CDFO}.
\newblock {\em European J. Oper. Res.}, 252(3):701--727, 2016.

\bibitem{bussieck2010minlp}
Michael~R Bussieck and Stefan Vigerske.
\newblock Minlp solver software.
\newblock {\em Wiley encyclopedia of operations research and management
  science}, 2010.

\bibitem{MR2312337}
Henri Cohen.
\newblock {\em Number theory. {V}ol. {I}. {T}ools and {D}iophantine equations},
  volume 239 of {\em Graduate Texts in Mathematics}.
\newblock Springer, New York, 2007.

\bibitem{MR2312338}
Henri Cohen.
\newblock {\em Number theory. {V}ol. {II}. {A}nalytic and modern tools}, volume
  240 of {\em Graduate Texts in Mathematics}.
\newblock Springer, New York, 2007.

\bibitem{MR3102373}
William Cook, Thorsten Koch, Daniel~E. Steffy, and Kati Wolter.
\newblock A hybrid branch-and-bound approach for exact rational mixed-integer
  programming.
\newblock {\em Math. Program. Comput.}, 5(3):305--344, 2013.

\bibitem{MR864015}
Wolfgang Dahmen.
\newblock Subdivision algorithms converge quadratically.
\newblock {\em J. Comput. Appl. Math.}, 16(2):145--158, 1986.

\bibitem{MR2869850}
Claudia D'Ambrosio and Andrea Lodi.
\newblock Mixed integer nonlinear programming tools: a practical overview.
\newblock {\em 4OR}, 9(4):329--349, 2011.

\bibitem{MR936450}
Carl de~Boor.
\newblock {$B$}-form basics.
\newblock In {\em Geometric modeling}, pages 131--148. SIAM, Philadelphia, PA,
  1987.

\bibitem{Dean:2008:MSD:1327452.1327492}
Jeffrey Dean and Sanjay Ghemawat.
\newblock Mapreduce: Simplified data processing on large clusters.
\newblock {\em Commun. ACM}, 51(1):107--113, January 2008.

\bibitem{Elber:2001:GCS:376957.376958}
Gershon Elber and Myung-Soo Kim.
\newblock Geometric constraint solver using multivariate rational spline
  functions.
\newblock In {\em Proceedings of the Sixth ACM Symposium on Solid Modeling and
  Applications}, SMA '01, pages 1--10, New York, NY, USA, 2001. ACM.

\bibitem{fourer}
R.~Fourer.
\newblock Linear programming: Software survey.
\newblock {\em OR/MS TODAY}, 42(3), 2015.

\bibitem{garloff1993bernstein}
J{\"u}rgen Garloff.
\newblock The bernstein algorithm.
\newblock {\em Interval computations (now Reliable Computing)}, 2(6):154--168,
  1993.

\bibitem{MR2059709}
Geoffrey~R. Grimmett and David~R. Stirzaker.
\newblock {\em Probability and random processes}.
\newblock Oxford University Press, New York, third edition, 2001.

\bibitem{MR0163361}
Theodore~E. Harris.
\newblock {\em The theory of branching processes}.
\newblock Die Grundlehren der Mathematischen Wissenschaften, Bd. 119.
  Springer-Verlag, Berlin; Prentice-Hall, Inc., Englewood Cliffs, N.J., 1963.

\bibitem{jeroslow}
R.~C. Jeroslow.
\newblock There cannot be any algorithm for integer programming with quadratic
  constraints.
\newblock {\em Operations Research}, 21(1):221--224, 1973.

\bibitem{book:747990}
Sven Leyffer~(eds.) Jon~Lee.
\newblock {\em Mixed Integer Nonlinear Programming}.
\newblock The IMA Volumes in Mathematics and its Applications 154.
  Springer-Verlag New York, 1 edition, 2012.

\bibitem{MR2640549}
Michael J{\"u}nger, Thomas Liebling, Denis Naddef, George Nemhauser, William
  Pulleyblank, Gerhard Reinelt, Giovanni Rinaldi, and Laurence Wolsey, editors.
\newblock {\em 50 years of integer programming 1958--2008}.
\newblock Springer-Verlag, Berlin, 2010.
\newblock From the early years to the state-of-the-art, Papers from the 12th
  Combinatorial Optimization Workshop (AUSSOIS 2008) held in Aussois, January
  7--11, 2008.

\bibitem{MR1244324}
Yuri~V. Matiyasevich.
\newblock {\em Hilbert's tenth problem}.
\newblock Foundations of Computing Series. MIT Press, Cambridge, MA, 1993.
\newblock Translated from the 1993 Russian original by the author, With a
  foreword by Martin Davis.

\bibitem{MR2499511}
B.~Mourrain and J.~P. Pavone.
\newblock Subdivision methods for solving polynomial equations.
\newblock {\em J. Symbolic Comput.}, 44(3):292--306, 2009.

\bibitem{MR3079719}
C{\'e}sar Mu{\~n}oz and Anthony Narkawicz.
\newblock Formalization of {B}ernstein polynomials and applications to global
  optimization.
\newblock {\em J. Automat. Reason.}, 51(2):151--196, 2013.

\bibitem{nataraj2007new}
P~SV Nataraj and M~Arounassalame.
\newblock A new subdivision algorithm for the bernstein polynomial approach to
  global optimization.
\newblock {\em International journal of automation and computing},
  4(4):342--352, 2007.

\bibitem{patil2014improved}
Bhagyesh~V Patil and PSV Nataraj.
\newblock An improved bernstein global optimization algorithm for minlp
  problems with application in process industry.
\newblock {\em Mathematics in Computer Science}, 8(3-4):357--377, 2014.

\bibitem{MR1266028}
Hartmut Prautzsch and Leif Kobbelt.
\newblock Convergence of subdivision and degree elevation.
\newblock {\em Adv. Comput. Math.}, 2(1):143--154, 1994.

\bibitem{Rajaraman:2011:MMD:2124405}
Anand Rajaraman and Jeffrey~David Ullman.
\newblock {\em Mining of Massive Datasets}.
\newblock Cambridge University Press, New York, NY, USA, 2011.

\bibitem{MR1263591}
Fabio Schoen.
\newblock Stochastic techniques for global optimization: a survey of recent
  advances.
\newblock {\em J. Global Optim.}, 1(3):207--228, 1991.

\bibitem{MR1689189}
Nigel~P. Smart.
\newblock {\em The algorithmic resolution of {D}iophantine equations},
  volume~41 of {\em London Mathematical Society Student Texts}.
\newblock Cambridge University Press, Cambridge, 1998.

\bibitem{MR1668773}
E.~M.~B. Smith and C.~C. Pantelides.
\newblock Symbolic methods for global optimisation.
\newblock In {\em Symbolic methods in control system analysis and design},
  volume~56 of {\em IEE Control Eng. Ser.}, pages 321--338. IEE, London, 1999.

\bibitem{MR0027121}
M.~H. Stone.
\newblock The generalized {W}eierstrass approximation theorem.
\newblock {\em Math. Mag.}, 21:167--184, 237--254, 1948.

\bibitem{sage}
{The Sage Developers}.
\newblock {\em {S}ageMath, the {S}age {M}athematics {S}oftware {S}ystem
  ({V}ersion 7.0)}, 2016.
\newblock {\tt http://www.sagemath.org}.

\bibitem{MR2496414}
Steve Zelditch.
\newblock Bernstein polynomials, {B}ergman kernels and toric {K}\"ahler
  varieties.
\newblock {\em J. Symplectic Geom.}, 7(2):51--76, 2009.

\end{thebibliography}

\end{document}